 \documentclass[11pt]{article}
 \usepackage[top=1in,bottom=1in,left=1.5in,right=1.5in]{geometry}
  \usepackage{amsmath,amssymb}
  \usepackage{latexsym}
 \usepackage[dvips]{pstricks} 
  \usepackage{pst-node}
  \usepackage[dvips]{graphicx}

  \newcommand{\C}{\mathbb{C}}
  
  \newcommand{\N}{\mathbb{N}}
  
  \newcommand{\R}{\mathbb{R}}
  \newcommand{\Z}{\mathbb{Z}}

  \renewcommand{\u}{\mathbf{u}}

  \newcommand{\x}{\mathbf{x}}

  \newcommand{\z}{\mathbf{z}}
  \newcommand{\0}{\mathbf{0}}
  \newcommand{\1}{\mathbf{1}}

  \newcommand{\cG}{\mathcal{G}}
  
  \newcommand{\cM}{\mathcal{M}}
  
  \newcommand{\cP}{\mathcal{P}}

  \newcommand{\Cp}{\mathrm{Cap\;}}

  \def\diag{\mathop{{\rm diag}}\nolimits}
  \newcommand{\hs}{\hspace*{\parindent}}

  \newcommand{\trans}{^\top}
  
  \newcommand{\per}{\mathop{\mathrm{perm}}\nolimits}
  \newcommand{\haff}{\mathrm{haf}}

  \newcommand{\rS}{\mathrm{S}}

  \newcommand{\perm}{\mathrm{perm}}

  \newcommand{\upp}{\mathop{\mathrm{upp}}\nolimits}
  
  \newtheorem{theo}{\bfseries \hs Theorem}[section]

  \newtheorem{lemma}[theo]{\bfseries \hs Lemma}
  
  \newtheorem{corol}[theo]{\bfseries \hs Corollary}
  \newtheorem{con}[theo]{\bfseries \hs Conjecture}

  \numberwithin{equation}{section} 

\begin{document}

 \title{Results and open problems in\\matchings in regular graphs\footnotemark[1]}
 \footnotetext[1]{This is an expanded version of the lecture ''Some open problems in matchings in graphs",
 Directions in Matrix Theory 2011, Coimbra, July 9, 2011, and AIM workshop ''Stability, hyperbolicity, and zero localization of functions",
 December 5--9, 2011.  The author acknowledges the support of Department of Mathematics, University of Coimbra
 and the Centre for Mathematics, University of Coimbra  and the American Institute of Mathematics, Palo Alto.}

 \author
 {Shmuel Friedland\\
 Department of Mathematics, Statistics and Computer Science\\
 University of Illinois at Chicago\\
 Chicago, Illinois 60607-7045, USA\\
 \emph{email}:friedlan@uic.edu}

 \date{3 January, 2012}
 \maketitle

 \begin{abstract}
 This survey paper deals with upper and lower bounds on the number of $k$-matchings in regular graphs on $N$ vertices.
 For the upper bounds we recall the upper matching conjecture which is known to hold for perfect matchings.
 For the lower bounds we first survey the known results for bipartite graphs, and their continuous versions as the van der Waerden and Tverberg permanent conjectures and its variants.  We then discuss non-bipartite graphs.  Little is known beyond the
 recent proof of the Lov\'asz-Plummer  conjecture on the exponential growth of perfect matchings  in cubic bridgeless graphs.
 We discuss the problem of the minimum of haffnians on the convex set
 of matrices, whose extreme points are the adjacency matrices of subgraphs of the complete graph corresponding to perfect matchings.
 We also consider infinite regular graphs.  The analog of $k$-matching is the $p$-monomer entropy, where $p\in [0,1]$ is the density
 of the number of matchings.

 \end{abstract}



 \section{Introduction}  Let $G=(V,E)$ be an undirected graph with no loops but possibly with multi-edges.
 We will identify $V$ with $[N]:=\{1,2,\ldots,N\}$, where $N=|V|$.
 A subset of edges $M\subset E$ is called a \emph{match} if no two edges in $M$ have a common vertex.  In physics literature
 $e:=(i,j)\in E$ is  called a \emph{dimer}.  $M$ is called a $k$-match if $\#M$, the number of edges in $M$, is $k$.
 Denote by $\cM_k$ the set of all $k$-matches of $G$.  Let $M\in\cM_k$.  So $M$ consists of $k$-dimers, and $M$ \emph{covers}
 $2k$ vertices, which the vertices of the induced graph $G(M)$ by $M$.
 $M$ is called a monomer-dimer cover of $G$, i.e. it consists of dimers which cover $2k$ vertices, while other vertices are called monomers.
 $M$ is called a \emph{perfect} matching if $M$ covers all vertices of $G$.  (So $V$ has to be even for $G$ to have a perfect matching.)
 We will usually assume that $|V|=2n$ unless stated otherwise.
 Let $\phi(k,G):=\#\cM_k(G)$ be the number of $k$-matchings in $G$.  Define $\phi(0,G):=1$.

 Denote by $\cG(2n,r)\supset\cG_{co}(2n,r)\supset\cG_{bi}(2n,r)$ the set of simple $r$-regular graphs on $2n$ vertices, the subset
 of $r$-colorable graphs, and the subset of bipartite $r$-regular graphs respectively.  Similarly denote by
 $\cG_{mult}(2n,r)\supset\cG_{co,mult}(2n,r)\supset\cG_{bi,mult}(2n,r)$  the set of $r$-regular graphs on $2n$ vertices, the subset
 of $r$-colorable graphs, and the subset of bipartite $r$-regular graphs respectively, where we allow multiple edges.

 The purpose of this survey paper is to report mostly on the following maximal and minimal problems
 \begin{eqnarray}\label{maxprobnb}
 \Theta(k,n,r):=\max\{\phi(k,G),\;G\in \cG(2n,r)\},\\
 \label{maxprobb}
 \Lambda(k,n,r):=\max\{\phi(k,G),\;G\in \cG_{bi}(2n,r)\},\\
 \label{minprobnb}
 \theta(k,n,r):=\min\{\phi(k,G),\;G\in \cG_{mult}(2n,r)\},\\
 \label{minprobnbco}
 \omega(k,n,r):=\min\{\phi(k,G),\;G\in \cG_{co,mult}(2n,r)\},\\
 \label{minprobb}
 \lambda(k,n,r):=\min\{\phi(k,G),\;G\in \cG_{bi,mult}(2n,r)\}.
 \end{eqnarray}

 It is known that for any simple graph $G=(V,E)$, with $|V|=2n$ an no isolated vertices, one has the following upper bound for the
 number of perfect matchings \cite{AF08, Ego07}
 \begin{equation}\label{AEFin}
 \phi(n,G)\le \prod_{v\in V}((\deg v)!)^{\frac{1}{2\deg(v)}},
 \end{equation}
 where $deg(v)$ is the degree of the vertex $V$.  Equality holds if and only if $G$ is a disjoint union of complete bipartite graphs.
 For a graph $G$ and $q\in\N$ denote by $qG$ the disjoint union of $q$ copies of $G$.  Let $K_{r,r}\in\cG_{bi}(2r,r)$ be a complete bipartite
 graph on $2k$ vertices.
 \eqref{AEFin} yields that $\Theta(qr,qr,r)=\phi(qr,qK_{r,r})=(r!)^q$.  The \emph{upper matching conjecture}, abbreviated as UMC, states
 \cite{FKM}
 \begin{equation}\label{UMC}
 \Theta(k,qr,r)=\phi(k,qK_{r,r}) \textrm{ for } k=1,\ldots,qr.
 \end{equation}
 This conjecture holds in the following cases: $k=1,2,3,4,qr$, for $r=2$ and for $r=3$ and $q=1,2,3,4$ \cite{FKM}.

 We know a lot about the lower bounds on the number of $k$-matchings in regular bipartite graphs.
 The main reason is the validity of the van der Waerden and Tverberg permanent conjectures \cite{vdW,Tve,Fr79,Ego, Fal}.
 The \emph{lower matching conjecture}, abbreviated here as LMC, states \cite{FKM}
 \begin{equation}\label{mmatchincon}
 \lambda(k,n,r)\ge {n\choose k}^2 (\frac{nr-k}{nr})^{rn-k}
 (\frac{kr}{n})^k \textrm{ for } k=1,\ldots,n.
 \end{equation}
 For $k=n$ the above conjecture follows from the lower bounds of \cite{Vor} for $r=3$ and from \cite{Sch} for any $r$.
 Weaker type of lower bounds can be found in \cite{FP,FG08,Gur11}.  A stronger version of the LMC is stated in \cite[(7.1)]{FKM}.
 For $r=2$ it is shown in \cite{FKM} that
 \begin{equation}\label{forlamkr=2}
 \lambda(k,n,2)=\phi(k,C_{2n}) \textrm{ for } k=1,\ldots,n,
 \end{equation}
 where $C_m$ is a cycle of length $m$.

 We do not know much about $\theta(k,n,r)$ and $\omega(k,n,r)$, except the case $r=2$ \cite{FKM}.
 It is straightforward to see that $\cG_{co,mult}(2n,2)=\cG_{bi,mult}(2n,2)$.  Hence $\omega(k,n,2)=\lambda(k,n,2)$ for $k=1,\ldots,n$.
 Furthermore \cite{FKM}
 \begin{eqnarray}\label{forthetkr=2a}
 \theta(k,n,2)=\phi(k,q C_3) \textrm{ for } k=1,\ldots,n, \textrm{ if } 2n=3q,\\
 \label{forthetkr=2b}
 \theta(k,n,2)=\phi(k,q C_3\cup C_4) \textrm{ for } k=1,\ldots,n, \textrm{ if } 2n=3q+4,\\
 \label{forthetkr=2c}
 \theta(k,n,2)=\phi(k,q C_3\cup C_5) \textrm{ for } k=1,\ldots,n, \textrm{ if } 2n=3q+5
 \end{eqnarray}

 It is well known that non-bipartite $r$-regular graphs on an even number
 of vertices may not have a perfect matching for $r\ge 2$.  For example, a simple $2$-regular graph that have exactly even number of odd cycles
 does not have a perfect matching.  Hence there exists $2r$-regular graphs with multiple edges with $2n$ vertices, for $n\ge 4$, which do not have
 perfect matches.  We will show that for any graph $G\in \cG_{mult}(2n,r)$ $\phi(k,G)\ge 1$ for $1\le k\le \frac{2n}{3}$.

 Assume that $G$ is a bridgeless $3$-regular graph on $2n$ vertices.
 Petersen's theorem implies that $G$ has a perfect matching.
  The Lov\'asz-Plummer conjecture states that
 there exists $\varepsilon>0$ such that \cite[\S8.7]{LP}
 \begin{equation}\label{LPcon}
 \phi(n,G) \ge 2^{\varepsilon 2n}.
 \end{equation}
 The Lov\'asz-Plummer conjecture was proved in \cite{EKKKN} with $\varepsilon=\frac{1}{3656}$.
 \begin{con}
 $\omega(n,n,r)$ has an exponential growth in $n$ for a fixed $r\ge 3$.
 \end{con}

 For an infinite graph $G=(V,E)$, with a countable number of vertices, match $M\subset E$ is
 a match of density $p\in [0,1]$, if  the proportion of vertices
 in $V$ covered by $M$ is $p$.
 Then the $p$-matching entropy of $G$ is defined as
 $$h_G(p)=
 \limsup_{k\to\infty}\frac{\log\phi(m_k,G_k)}{\#V_k},$$
 where $G_k=(E_k,V_k), k\in\N$ is a sequence of finite graphs
 converging to $G$, and $\lim_{k\to\infty}
 \frac{2m_k}{\#V_k}=p$.  See for details \cite{FKLM}.
 Denote by $\infty K_{r,r}$ a countable union of disjoint copies of $K_{r,r}$.
 One can compute $h_{\infty K_{r,r}}(p)$ for any $p\in[0,1]$ using the thermodynamics formalism \cite{FKLM}.
 Assume that $G$ is an infinite $r$-regular graph.
 The \emph{asymptotic matching conjecture}, abbreviated as AUMC, states that $G$ satisfies
 \begin{equation}\label{AUMC}
 h_G(p)\le h_{\infty K_{r,r}}(p) \textrm{ for all } p\in[0,1].
 \end{equation}
 For bipartite $r$-regular graphs infinite graphs $G$ the AUMC was stated in \cite{FKLM}.
 Since the UMC holds for perfect matchings it follows that \eqref{AUMC} holds for $p=1$.  (See for details the arguments in \cite{FKLM}.)

 Let
 \begin{equation}\label{deffrp}
  f(p,r):=\frac{1}{2}\big(p \log r -p\log p - 2(1-p)\log (1-p)
 +(r-p)\log (1 -\frac{p}{r})\big).
 \end{equation}
 The \emph{asymptotic lower matching conjecture}, abbreviated ALMC, claims that for any $r$-regular bipartite infinite graph $G$ the following
 inequality holds  \cite{FKM}.
 \begin{equation}\label{almc}
 h_G(p)\ge f(p,r) \textrm{ for all } p\in [0,1].
 \end{equation}
 This inequality follows straightforward from the LMC.
 For $r=2$ holds this conjecture holds in view of \eqref{forlamkr=2}.  Equality holds if $G$ is an infinite path, i.e. the monomer-dimer lattice
 on $\Z$.  For $r\ge 3$ the conjecture holds for all $p$ of the form $\frac{r}{r+s}, s=0,1,\ldots,$ \cite{FG08}.
 A weaker version of the ALMC is proved in \cite{FKLM}.
 For $r=2d$ the ALMC holds for the monomer-dimer model on $\Z^d$ for $d=2,3,\ldots,$ \cite{FF11}.  The complete proof of the ALMC is given in
 \cite{Gur11}.
 There should be a similar inequality for infinite $r$-regular graphs, at least for $p\in [0,\frac{2}{3}]$ as we explain in \S7.
 However we do not have the precise version of such inequality.

 We now survey briefly the contents of this paper.  In \S2 we discuss the notions of matching polynomials, haffnians and permanents.
 In \S3 we give the known upper bounds for matchings in regular graphs.  In \S4 we discuss the proved van der Waerden and Tverberg
 conjectures, their discreet analogs and their relation to the lower estimates on the number of $k$-matchings in multi $r$-regular bipartite
 graphs.  In \S5 we give a brief introduction to positive hyperbolic polynomials which yield short elegant proofs to lower bounds discussed
 in \S4.  In \S6 we discuss matching in non-bipartite graphs.  Here we mostly have open problems.  In \S7 we bring the known results
 and conjectures on the entropy of bipartite infinite $r$-regular, and similar open problems for infinite non-bipartite $r$-regular graphs.

 \section{Preliminary results}
  Denote by $\rS_0(N,\R_+)$ the set of all $N\times N$ symmetric matrices $W=[w_{ij}]_{i=j=1}^N$ with zero diagonal and nonnegative entries.
 Let $G=(V,E)$ be a simple undirected graph, where $N=|V|$ is an arbitrary positive integer.  $W(G)=[w_{ij}]_{i=j=1}^N\in\rS_0(N,\R_+)$ is called a weighted adjacency matrix if $ w_{ij}=0$ if $(i,j)\not\in E$.  The adjacency matrix $A(G)$ is the maximal $0-1$ matrix which a weighted adjacency matrix of $G$.  Note that any $W\in \rS_0(N,\R_+)$ is a weighted adjacency matrix of $K_N$, the complete graph on $N$ vertices.
 The multivariate matching polynomial of $K_N$ is a multilinear polynomial in $N$ complex variables $\z=(z_1,\ldots,z_N)\in\C^N$:
 \begin{equation}\label{marchpolNvar}
 \Psi(\z,W):=1+\sum_{k=1}^{\lfloor\frac{N}{2}\rfloor} \sum_{M\in\cM_k(K_N)}\prod_{(i,j)\in M}w_{ij}z_iz_j
 \end{equation}
 The fundamental result of Heilmann and Lieb \cite{HL} claims that $\Psi(\z,W)$ is nonzero if either $\Re(z_i)>0, i=1,\ldots,N$ or $\Re(z_i)<0, i=1,\ldots,N$. (We call this property the left-right half plane property.)  Another proof of this fact is given in \cite{COSW}, where it is observed that the following polynomial has the left-right half plane property:
 \begin{equation}\label{sokpol}
 \Upsilon(\z,W)=\prod_{1\le i < j\le N} (1+w_{ij}z_iz_j).
 \end{equation}
 From this fact it is deduced in \cite{COSW} that $\Psi(\z,W)$ has the left-right half plane property.

 Recall that the $k$-haffnian of $W$ is given by
 \begin{equation}\label{defkhaf}
 \haff_k W:=\sum_{M\in\cM_k(K_N)} \prod_{(i,j)\in M}w_{ij}.
 \end{equation}
 Then a matching polynomial of $\Phi(t)$ is given by
 \begin{equation}\label{defmatpol}
 \Phi(t,W))=1+\sum_{k=1}^{\lfloor\frac{N}{2}\rfloor} \haff_k W t^k.
 \end{equation}
 Clearly $\Psi((z,\ldots,z),W)=\Phi(z^2,W)$.  The left-right hyperplane property yields that $\Phi(z^2,W)=0$ implies that $\Re(z)=0$.
 Hence $\Psi(t,W)$ has only real negative roots.  Apply the Newton inequality to deduce
 \begin{equation}\label{newtin}
 (\frac{\haff_l W }{{\lfloor\frac{N}{2}\rfloor\choose l}})^2 \ge \frac{\haff_p W }{{\lfloor\frac{N}{2}\rfloor\choose p}}
 \frac{\haff_q W}{{\lfloor\frac{N}{2}\rfloor\choose q}} \textrm{ for }1\le p <l<q \le \lfloor\frac{N}{2}\rfloor.
 \end{equation}

 Let $W=A(G)$ is the adjacency matrix of a simple graph on $G$ on $N$ vertices.  Then $\haff_k A(G)$ is the number of $k$-matchings in $G$.  Denote by $\rS_0(N,\Z_+)$ the set of symmetric matrices with nonnegative integer entries and zero diagonal.
 Then $A=[a_{ij}]\in\rS(N,\Z_+)$ is the adjacency matrix of a multigraph $G(A)$. i.e. $a_{ij}$ is the multiplicity of the edge $(i,j)$.
 So $\haff_k A$ is the number of $k$-matches in $G(A)$.

 Let $K_{m,n}$ be a complete bipartite graph $G=(V_1\cup V_2,E)$ where $V_1=[m],V_2=[n]$.  Then
 \[A(K_{m,n})=\left[\begin{array}{cc}0_{m\times n}&1_{m \times n}\\ 1_{n\times m}&0_{n \times m}\end{array}\right],\]
 where $1_{m\times n}$ is an $m\times n$ matrix whose entries are $1$.  Then $W(K_{m,n})$, the weighted adjacency matrix of $K_{m,n}$
 has the form
 \[W=\left[\begin{array}{cc}0_{m\times n}&B\\ B\trans &0_{n \times m}\end{array}\right], \quad B\in\R^{m\times n}_+.\]
 $B$ is called the weighted bipartite adjacency matrix of $K_{m,n}$.  Similary, for a bipartite graph $G=([m]\cup [n],E)$,
 the matrix $B=[b_{ij}]\in \R^{m\times n}_+$ is called the weighted bipartite adjacency matrix of $G$ if $b_{ij}=0$ for $i\in [m],j\in [n]$.
 Assume that $m=n$.  Then $\haff_n W$ is equal to the permanent of $B$, denoted by $\perm\; B$ or $\perm_n B$.
 For any $B\in\R_+^{m\times n}$ and a positive integer $k\le min(m,n)$ denote by $\perm_k B$ the sum of the permanents of all $k\times k$ submatrices of $B$.  Then $\haff_k W =\per_k B$.

 \section{Upper bounds on matchings}
 Denote by $\Z_+^{n\times n}\supset \{0,1\}^{n\times n}$ the set of all $n\times n$ matrices with nonnegative integer entries,
 and the subset of all $0-1$ matrices.  Each $B\in\Z_+^{n\times n}$ represents a bipartite multigraph $G=([n]\cup[n],E)$, where $B$ is a bipartite
 adjacency matrix of $G$.  $G$ is simple if and only if $B\in\{0,1\}^{n\times n}$.  Let $r_i(B)$ the $i$-th row sum of $B$, i.e. the degree
 of the vertex $i\in V_1=[n]$.   Bregman's inequality \cite{Bre}, conjectured by Minc \cite{Min63}, states
 \begin{equation}\label{bregin}
 \per B\le \prod_{i=1}^n (r_i(B)!)^{\frac{1}{r_i(B)}}, \quad B\in \{0,1\}^{n\times n}.
 \end{equation}
 Equality holds if and only if $B$ the bipartite adjacency matrix of union of complete bipartite graphs.
 Let $G=([2n],E)$ be a simple graph, with the adjacency matrix $A(G)\in\rS_0(2n,\{0,1\})\subset \{0,1\}^{n\times n}$.
 Consider the bipartite graph $G'=(V_1\cup V_2, E')$ on $4n$ vertices whose bipartite adjacency matrix is $A(G)$.
 The square of the number of perfect matchings of $G$ counts
 ordered pairs of such matchings. We claim that this is the number of spanning $2$-regular 
 subgraphs $H$ of $G$ consisting of even cycles (including cycles of length $2$
 which are the same edge taken twice), where each such $H$ is counted $2^s$ times, with
 $s$ being the number of components (that is, cycles) of H with more than $2$ vertices.
 Indeed, every union of a pair of perfect matchings $M_1,M_2$ is a $2$-regular spanning
 subgraph $H$ as above, and for every cycle of length exceeding $2$ in H there are two
 ways to decide which edges came from $M_1$ and which from $M_2$.
 The permanent of $A(G)$ also counts the number of spanning
 $2$-regular subgraphs $H'$ of $G$, where now we allow odd cycles and cycles of
 length $2$ as well. Here, too, each such $H'$ is counted $2^s$ times, where $s$ is the number
 of cycles of $H'$ with more than $2$ vertices, (as there are $2$ ways to orient each such
 cycle as a directed cycle and get a contribution to the permanent). Thus the square
 of the number of perfect matchings is at most the permanent of $A(G)$.
 Use \eqref{bregin} to deduce \eqref{AEFin}.  Equality holds if and only if $G$ is a disjoint union of complete bipartite
 graphs \cite{AF08}.  Hence
 \begin{equation}\label{upperfmatch}
 \phi(n,G)\le (r!)^{\frac{n}{r}} \textrm{ for }G\in \cG(2n,r).
 \end{equation}
 Equality holds if $n=qr$ and $G=qK_{r,r}$.
 The above inequality imply the following inequalities:
 \begin{eqnarray}\label{uppkmatch}
 \Theta(k,n,r)\le \min({2n \choose 2k}(r!)^{\frac{k}{r}},2^{-k}{2n \choose k} r^k) \textrm{ for }  k=1,\ldots,n\\
 \label{uppkmatchbi}
 \Lambda(k,n,r)\le \min({n \choose k}^2(r!)^{\frac{k}{r}},{n \choose k} r^k) \textrm{ for } k=1,\ldots,n.
 \end{eqnarray}
 Since the inequalities \eqref{uppkmatchbi} are shown in \cite{FKM}, we give only a simple proof of \eqref{uppkmatch}.
 Let $G \in\cG(2n,r)$.  Fix $k\in [1,n]$ and choose a subgraph $H$ of $G$ on $2k$ vertices.  Use \eqref{AEFin} to deduce that
 $\phi(k,H)\le (r!)^{\frac{k}{r}}$.  Hence $\phi(k,G)\le {2n \choose 2k}(r!)^{\frac{k}{r}}$.
 The second inequality is obtained as follows.  Choose $k$ vertices $v_1,\ldots,v_k$ in $G$ and consider the maximum number of $k$ matches
 of the form $(v_1,w_1),\ldots,(v_k,w_k)$.  Since the degree of each $v_i$ is $r$ the number of such $k$-matches is at most $r^k$.
 The number of choices of $k$ vertices in $G$ is $2n \choose k$.  Note that a $k$-match of the form $(v_1,w_1),\ldots,(v_k,w_k)$
 is counted exactly $2^k$ times.  Hence $\phi(k,G)\le 2^{-k}{2n \choose k}r^k$.

 \section{Van der Waerden and Tverberg conjectures and their discreet analogs}
 Let $\Omega_n\subset \R^{n\times n}$ be the convex set of doubly stochastic matrices.  Birkhoff-von Neumann theorem claims that the set of the extreme points of $\Omega_n$ is the set of permutation matrices $\cP_n\subset \Omega_n$.   Note that if $G\in \cG_{bi,mult}(2n,r)$ then $\frac{1}{r}B(G)\in\Omega_n$.  Birkhoff theorem yields that $G(B)$ is a sum of $r$-permutation matrices for each $G\in \cG_{bi,mult}(2n,r)$.  So each $G\in \cG_{bi,mult}(2n,r)$ is $r$-colorable.  Thus
 \begin{equation}\label{basbimatin}
 \min\{\perm_k(C), C\in\Omega_n\}\le\frac{\lambda(k,n,r)}{r^k} \textrm{ for } k=1,\ldots,n.
 \end{equation}
 It is known that
 \begin{equation}\label{vwdtvercon}
 \mu(k,n):=\min\{\perm_k(C), C\in\Omega_n\}=\perm_k (\frac{1}{n} 1_{n\times n})={n\choose k}^2 \frac{k!}{n^r} \textrm{ for }  k=1,\ldots, n.
 \end{equation}
 For $k=n$ this equality was the celebrated van der Waerden conjecture \cite{vdW}.  An inequality $\mu(n,n)\ge e^{-n}$ was given in \cite{Fr79}.
 This inequality proved the Erd\H{o}s-R\'enyi conjecture that the number of matchings in $3$-regular bipartite graphs grows exponentially.
 The van der Waerden conjecture was proved by Egorichev \cite{Ego} and Falikman \cite{Fal}.  The Tverberg conjecture was proved by the author
 in \cite{Fr82}.

 It was shown by Schrijver \cite{Sch} that
 \begin{equation}\label{schirjin}
 \lambda(n,n,r)\ge (\frac{(r-1)^{r-1}}{r^{r-2}})^n.
 \end{equation}
 The case $r=3$ was shown by Voorhoeve \cite{Vor}, which also gave another proof to the above Erd\"os-R\'enyi conjecture.
 Gurvits \cite{Gur06} gave an improved continuous version of the above inequality.  Assume that
 $B\in\Omega_n$ has in each column has at most $r$ nonzero entries. Then
 \begin{equation}\label{Gurvin}
 \per B \ge \frac{r!}{r^r}\big(\frac{r}{r-1}
 \big)^{r(r-1)}\big(\frac{r-1}{r}\big)^{(r-1)n}.
 \end{equation}
 Hence
 \begin{equation}\label{Gurvin1}
 \lambda(n,n,r)\ge  \frac{r!}{r^r}\big(\frac{r}{r-1}
 \big)^{r(r-1)}(\frac{(r-1)^{r-1}}{r^{r-2}})^n.
 \end{equation}
 See \cite{LS} for a simple exposition of Gurvits's result.
 In the joint work with Gurvits we give the following lower bound \cite[Theorem 4.5]{FG08}
 \begin{equation}\label{lbfg08}
 \lambda(k,n,r)\ge \max_{s\in [n-r]} \frac{(s n)!{n\choose k}}{s^{n-k}(n-k)!((s-1)n+k)!}\frac{r^k(r+s)!}{(r+s)^{r+s}}
 \big(\frac{r+s-1}{r+s})^{(r+s-1)(n-r-s)}.
 \end{equation}

 \section{Positive hyperbolic polynomials and lower bounds on permanents}
 In this section we review briefly the main ideas and techniques to show \eqref{vwdtvercon}, \eqref{Gurvin1} and
 \eqref{lbfg08}.  This exposition relies on the papers \cite{FG06,FG08}.

 A polynomial $p=p(\x)=p(x_1,\ldots,x_n):\R^n
 \to\R$ is called \emph{positive hyperbolic} if the following conditions hold
 \begin{enumerate}
 \item
 $p$ is a homogeneous polynomial of degree $m\ge 0$.
 \item
 $p(\x)>0$ for all $\x>0$.
 \item
 $\phi(t):=p(\x+t\u)$, for $t\in \R$, has $m$-real $t$-roots for each
 $\u>\0$ and each $\x$.
 \end{enumerate}
 Assume that $p$ is positive hyperbolic.
 Then the coefficient of each monomial in positive hyperbolic polynomial is nonnegative.
 Moreover $p(x_1,\ldots,x_{n-1},0)$ and $\frac{\partial }{\partial x_n}p(x_1,\ldots,x_{n-1},0)$ are either
 positive hyperbolic or zero polynomials.  The following examples are most useful for positive hyperbolic polynomials.\medskip

 \noindent
 \emph{Example 1}:  Let $A=[a_{ij}]_{i=j=1}^{m,n}\in \R_+^{m\times n}$, where each row of $A$ is nonzero.
 Then $p_{k,A}(\x):=\sum_{1\le i_1<\ldots < i_k\le m} \prod_{j=1}^k(A\x)_{i_j},
 \x\in\R^n$ is positive hyperbolic.\medskip

 \noindent
 \emph{Example 2}:  Let $A_1,\ldots,A_n\in \C^{m\times m}$ be hermitian,
 nonnegative definite matrices such that $A_1+\ldots+A_n$ is
 a positive definite matrix.  Let $p(\x)=\det \sum_{i=1}^n x_i
 A_i$.  Then $p(\x)$ is positive hyperbolic.\medskip

 \noindent
 \emph{Example 3}: $B\in \R_+^{m\times m}$ symmetric.  Then $\x\trans B\x$ is positive hyperbolic if and only if
 $B$ has exactly one positive eigenvalue.\medskip

 Let $p(\x):\R^n\to \R$ positive hyperbolic polynomial of degree $m\ge 1$.
 Define the \emph{capacity} of $p$ as follows
 \begin{equation}\label{defcapp}
 \Cp p:=\inf_{\x> 0, x_1\cdots x_n =1} p(\x).
 \end{equation}
 Since $\log t$ is a concave function of $t$, use Birkhoff theorem for doubly stochastic matrices to deduce that
 \begin{equation}\label{capfords}
 \Cp p_{k,A}={n \choose k} \textrm{ for each }A\in\Omega_n.
 \end{equation}
 Let $B=D_1 A D_2$,
 $D_1,D_2$ positive definite diagonal matrices and $A$ doubly
 stochastic matrix.  Let $p_{n,B}$ be defined as above.
 Then $\Cp p_{n,B}=\frac{1}{\det D_1D_2}$.  Recall that any positive square matrix $B$ has the above decomposition \cite{Sin64,BPS66}.

 Let $A\in\R^{n\times n}_+$ and consider the hyperbolic polynomial $p(\x):=(x_1+\ldots+x_n)^{n-k}p_{k,A}$.  Then its mixed partial derivative
 is given by
 \begin{equation}\label{mixderfor}
 \frac{\partial ^n}{\partial x_1\ldots\partial x_n}p(\0)=(n-k)! \perm_k A,\quad p(\x):=(x_1+\ldots+x_n)^{n-k}p_{k,A}.
 \end{equation}

 This to estimate from below the mixed partial derivative of a positive hyperbolic polynomial $p$ we need to give a lower bound
 on the capacity of

 \noindent
 $\frac{\partial p}{\partial x_i}(x_1,\ldots,x_{i-1},0,x_{i+1},\ldots,x_n)$ in terms of the capacity of $p$.
 Denote by $\deg_i p$ the highest degree of $x_i$ in all nonzero monomials that appear in $p$.  The following lemma is the precise
 form of the lower bound explained above \cite[Lemma 2.4.]{FG08}.
 \begin{lemma}\label{bascapin} Let $p:\R^n\to \R$  be a positive hyperbolic
 of degree $m\ge 1$.  Assume that $\Cp p>0$. Then $\deg_i p \ge 1$ for $i=1,\ldots, n$.
 For $m=n\ge 2$
 $$\Cp \frac{\partial p}{\partial
 x_i}(x_1,\ldots,x_{i-1},0,x_{i+1},\ldots,x_n)\ge (\frac{\deg_i p
 -1}{\deg_i p})^{\deg_i p-1} \Cp p \textrm{ for }
 i=1,\ldots,n,$$
 where $0^0=1$.
 \end{lemma}

 Combine the above lemma with \eqref{capfords} to deduce the Tverberg conjecture \eqref{vwdtvercon}.
 The inequality \eqref{Gurvin} for any doubly stochastic matrix, which has in each column has at most $r$ nonzero entries,
 follows from the observation that $\deg_i p_{n,B}\le r$ for $i=1,\ldots,n$ \cite{Gur06}.

 The inequality \eqref{lbfg08} is obtained in the following way.  Assume $B\in\Omega_n$ in each column has at most $r$ nonzero entries.
 Fix a positive integer $s\in [n-k]$.  Choose a bipartite graph $G\in \cG_{bi,mult}(2n,s)$.  Let $B(G)$ be the bipartite adjacency matrix
 of $G$.  Now form the positive hyperbolic polynomial $p=p_{n-k,B(G)}p_{k,B}$.  Note that $\Cp p_{n-k,B(G)}p_{k,B}\ge {n \choose k}^2$, and
 $\deg_i p\le r+s$ for $i=1,\ldots,n$.  Hence the mixed derivative of $p$ is bounded below by the right-hand side of \eqref{Gurvin},
 where $r$ is replaced by $r+s$.  Now sum over all $G\in \cG_{bi,mult}(2n,s)$.  This will yield the inequality
 \begin{equation}\label{fginper}
 \perm_k B\ge \frac{(s n)!{n\choose k}}{s^{n-k}(n-k)!((s-1)n+k)!}\frac{(r+s)!}{(r+s)^{r+s}}
 \big(\frac{r+s-1}{r+s})^{(r+s-1)(n-r-s)},
 \end{equation}
 for $s=1,\ldots,n-k$.  This implies \eqref{lbfg08}.

  \section{Non-bipartite graphs}

 Let $\Omega_{n,s}\subset \Omega_n$ be the convex subset of all symmetric doubly stochastic matrices.
 The extreme points of $\Omega_{n,s}$ were determined by M. Katz \cite{Kat70}.
 \begin{theo}\label{katzchar}  $F\in\Omega_{n,s}$ is an extreme point of $\Omega_{n,s}$ if and only if the following condition holds.
 There exists a permutation matrix $P\in\Omega_n$ such that $PFP\trans$ is a block diagonal matrix  $\diag(F_1,\ldots,F_m)$ where each $F_i$ is one of the following three possible forms. (i) $F_i=[1]\in \Omega_1$; (ii) $F_i=\left[\begin{array}{cc}0&1\\1&0\end{array}\right]\in\Omega_2$;
 (iii) $F_i=\frac{1}{2}A(C_{2j-1})\in \Omega_{2j-1}$, where $C_{2j-1}$ is an odd cycle of length $2j-1$ for $j\ge 2$.
 \end{theo}
 \begin{corol}\label{katzcor} Let $\Omega_{n,s,0}$ be the convex set of symmetric doubly stochastic matrices with zero diagonal.
 Then each extreme point $F$ of $\Omega_{n,s,0}$ corresponds to a spanning graph subgraph $G$ of $K_n$ which is a union of $k$-match $M$
 and odd cycles.  $F$ is a weighted adjacency matrix of $G$, where the weight of each edge in $M$ is $1$, and the weight of each edge in an
 odd cycle is $\frac{1}{2}$.
 \end{corol}

 Since $C_{2j-1}$ has a $j$-match we deduce
 \begin{theo}\label{minmatcreggraph}  Let $G$ be $r$-multi-regular graph on $n$ vertices with no selfloops.
 Then $G$ has $\lceil\frac{n}{3}\rceil$-match.  Equality holds if $r=2l$, $n=3q$, and $G$ is a union of $q$ $2l$-multigraphs
 on three vertices without selfloops.
 In particular $\omega(k,n,r)>0$ for $k\in [\lceil\frac{2n}{3}\rceil]$.
 \end{theo}

 It is natural to conjecture that for any $p\in (0,\frac{2}{3})$ and an integer $k\in [pn,\frac{2n}{3}]$ $\omega(k,n,r)$ has an exponential
 growth in $n$ for a fixed $r\ge 2$.  (For $r=2$ this follows from \cite{FKM}.)  The real challenge is to give a good lower bound for the exponential growth.

 In the rest of this section we discuss the analogous problem to the van der Waerden and Tverberg permanent conjecture.
 Most of our exposition are from \cite{Fr11}.
 Let $\Psi_{2n}\subset \Omega_{2n,s,0}$ be the convex set of symmetric doubly stochastic matrices with
 zero diagonal, whose extreme points are the adjacency matrices of perfect matches in $K_{2n}$.  $\Psi_{2n}$ was characterized by Edmonds \cite{Ed65a}, see \cite[Theorem 6.3, 2nd Proof, page 209]{CCPS} for a simple proof.
 Namely, it is the set of all $B=[b_{ij}]\in \Omega_{2n,s,0}$ that satisfy the conditions
 \begin{equation}\label{edmcon}
 \sum_{i,j\in S} b_{ij}\le |S|-1, \textrm{ for each } S\subset{1,\ldots,2n}, |S| \textrm{ odd and } 3\le |S|\le 2n-3.
 \end{equation}

 Our problem is to find or give a good lower bound for
 \begin{equation}\label{fridprob}
 \min_{B\in\Psi_{2n}}\haff_k(B)= \mu_{k,n} \textrm{ for } k=2,\ldots,n.
 \end{equation}
 It is tempting to state, as in the case of the van der Waerden and Tverberg's conjectures that
 \begin{equation}\label{fridcon}
 \mu_{k,n}=\haff_k(\frac{1}{2n-1}A(K_{2n}))\textrm{ for } k=2,\ldots,n.
 \end{equation}
 Equality holds if and only if $B=\frac{1}{2n-1}A(K_{2n})$.
 (According to an e-mail from Leonid Gurvits, he stated this conjecture for $k=n$ in correspondence with E. Lieb on September 21, 2005.)
 It is easy to show that this conjecture is true for $k=2$.  However for $k=n$ and $n$ big enough \eqref{fridcon} is wrong as explained below.

 \cite[Lemma 2.1]{Fr11} states:
 \begin{lemma}\label{exactfor}  For positive integers $2\le k\le n$ we have
 \begin{eqnarray}\label{exactfor1h}
 \haff_k(\frac{1}{2n-1} A(K_{2n}))=\frac{1}{(2n-1)^k} {2n\choose 2k} \frac{1}{k!}\prod_{j=0}^{k-1} {2k-2j\choose 2},\\
 \haff_k(\frac{1}{n} A(K_{n,n}))=\frac{1}{n^k} {n \choose k}^2 k!\label{exactfor1p}.
 \end{eqnarray}
 In particular
 \begin{equation}\label{exactfor2}
  e^{-n}\sqrt{2}<\haff_n(\frac{1}{2n-1} A(K_{2n}))=\frac{(2n)!}{(2n-1)^n 2^n n!}<\haff_n(\frac{1}{n} A(K_{n,n}))=\frac{n!}{n^n}
 \end{equation}
 For $n\gg 1$ we have the following approximations
 \begin{equation}\label{approx}
 \haff_n(\frac{1}{2n-1} A(K_{2n}))\approx e^{-n}\sqrt{2e} , \quad \haff_n(\frac{1}{n} A(K_{n,n}))\approx e^{-n}\sqrt{2\pi n}.
 \end{equation}
 \end{lemma}

 Observe next that if $G=(V,E)$ is a $r$-regular graph without loops on an even number of vertices, then $\frac{1}{r}A(G)$ is in $\Psi_{|V|}$
 if and only if any vertex cut $S\subset V$ with an odd number of vertices, $3\le |S|\le |V|-3$, has at least $r$ edges.  Hence if Conjecture \eqref{fridcon} holds, then \eqref{approx} implies that such a regular graph has at least $(\frac{r}{e})^{\frac{|V|}{2}}$ perfect matchings.
 In \cite{CPS} the authors construct an infinite family of cubic 3-colored connected graph $G=(V,E)$, for which the number of matchings is less then $c_F|V|(\frac{1+\sqrt{5}}{2})^{\frac{|V|}{12}}$. (Here $|V|=12k+4$ and $k=1,2,\ldots$.)  As $(\frac{1+\sqrt{5}}{2})^{\frac{1}{12}}< 1.017< \sqrt{\frac{3}{e}}\approx 1.05$ we must have that $\mu_{n,n}<\haff (\frac{1}{2n-1}A(K_{2n}))$ for $n\gg 1$.  (I would like to thank S. Norin for pointing out to me this fact.)

 Since $\mu_{n,n}$ is the minimum of the haffnian function, it follows that $\mu_{n+m,n+m}\le \mu_{n,n}\mu_{m,m}$.  Hence the sequence $\log\mu_{n,n}$ is subadditive. In particular the following limit exists
 \begin{equation}\label{defmu}
 \mu:=\lim_{n\to\infty} \frac{\log\mu_{n,n}}{n}.
 \end{equation}
 A weak analog of the van der Waerden conjecture is the claim that $\mu>-\infty$.  Note that the above example in \cite{CPS} implies that $\mu\le \frac{\log\frac{1+\sqrt{5}}{2}}{6}-\log 3$.  Other generalizations of the van der Waerden conjectures for perfect matchings in hypergraphs are
 considered in \cite{BS}.

 We close this section with a lower bound of $\perm_n B$ for $B$ is a corresponding neighborhood of $A=\frac{1}{2n-1} A(K_{2n})$
 \cite[Theorem 2.2]{Fr11}.
 \begin{theo}\label{hypest}  Let $B\in \Psi_{2n}$. Assume that $B\in\Psi_{2n}$ has exactly one positive eigenvalue.  Then
 \begin{equation}\label{hypest1}
 \haff_n(B)\ge (\frac{n-1}{n})^{(n-1)n}\approx e^{-n}\sqrt{e}.
 \end{equation}
 Moreover, for each $k=2,\ldots,n-1$
 \begin{equation}\label{hypest2}
 \haff_k(B)\ge \frac{(2n)^{2n-2k} (2n-k)! (2n)^k}{(2n-2k)!(2n-k)^{2n-k}2^k k!} (\frac{(2n-k-1)}{2n-k})^{(2n-k-1)k}.
 \end{equation}

 \end{theo}
 The proof of the theorem follows from the fact that from Example 3 in \S5, which stated that $\x\trans B\x$ is positive hyperbolic polynomial.

 \section{Infinite regular graphs}
 Let $\alpha_{n,r}$ be a probability measure on $\cG_{bi,mult}(2n,r)$.
 There are two simple natural measures on $\cG_{bi,mult}(2n,r)$ \cite{FKM}.
 We describe one that easily generalizes to $\cG_{co,mult}(2n,r)$.
 For each $G\in\cG_{bi,mult}(2n,r)$ the adjacency matrix $A(G)$ is the sum of $r$-permutation matrices $P_1+\ldots+P_r$.
 Let $N(G)$ be the number of all possible representation of $A(G)$, viewed as $(P_1,\ldots,P_r)$.  Then
 $\alpha_{n,r}(G)=\frac{N(G)}{(n!)^r}$.
 Let $E(m,n,r)$ the expected number of $m$-matchings on $\cG_{bi,mult}(2n,r)$ with respect to $\alpha_{n,r}$.
 \cite[Theorem 4.4]{FKM} claims
 \begin{equation}\label{limval}
 \lim_{k\to\infty} \frac{\log E(m_k,n_k,r)}{2n_k}=\frac{1}{2}(p \log
 r -p\log p - 2(1-p)\log (1-p)
 +(r-p)\log (1 -\frac{p}{r})).
 \end{equation}
 We assume here that $1\le m_k \le n_k,
 k=1,...,$ are two strictly increasing sequences of integers such
 that the sequence $\frac{m_k}{n_k}, k=1,...$ converges to $p\in
 [0,1]$.    The above equality is the motivation for the ALMC \eqref{almc}.

  Figure 1 gives the graphs of the AUMC, ALMC, the $20$ points of $h_{\Z^2}(p)$, and the lower bound $FT$ given by the Tverberg permanent conjecture for $r=4$.  ($\Z^2$ is the infinite graph, whose vertices are the integer points on $\R^2$, and there is an edge between the points
  the integer points $(i,j), (p,q)$ if and only of $|i-j|+|p-q|=1$.)  The $20$ points of $h_{\Z^2}(p)$  we first computed by Baxter \cite{Bax}
  using his famous corner transfer matrix, and reconfirmed in \cite{FP11} using upper and lower bounds for $h_{\Z^2}(p)$.
  The value $h_2:=\max_{p\in[0,1]} h_{\Z^2}(p)$ is the $2$-dimensional monomer-dimer entropy of $\Z^2$ computed in \cite{FP}.
  \begin{figure}[here]
    \begin{center}
    \includegraphics[width=0.75\textwidth]{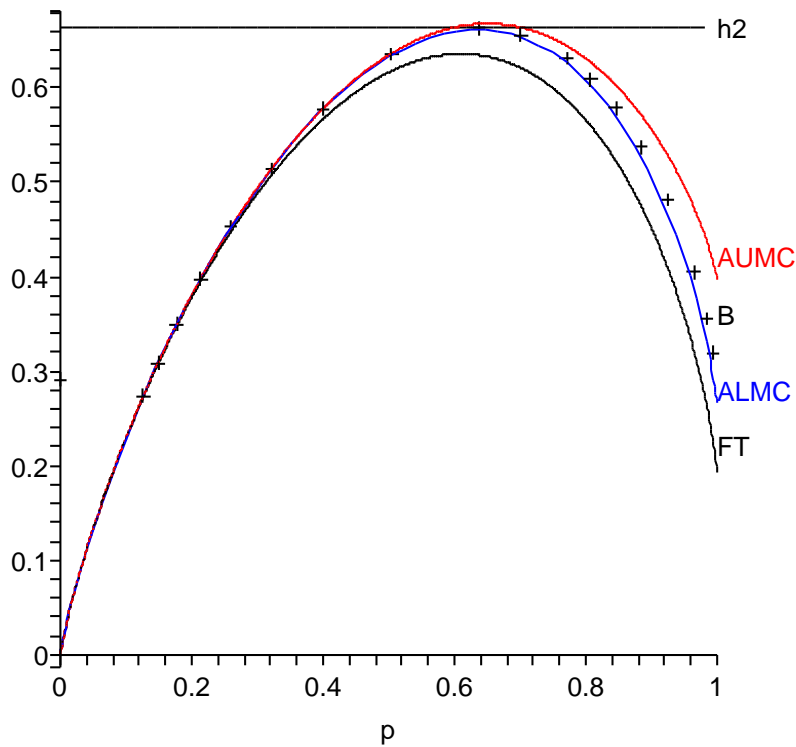}
    \caption{Monomer-dimer tiling of the $2$-dimensional grid:
    entropy as a function of dimer density. FT is the
    Friedland-Tverberg lower bound, h2 is the true monomer-dimer
    entropy. B are Baxter's computed values. ALMC is the Asymptotic
    Lower Matching Conjecture. AUMC is the entropy of a countable
    union of $K_{4,4}$, conjectured to be an upper bound by the
    Asymptotic Upper Matching Conjecture.}\label{Fig:h2bounds}
    \end{center}
  \end{figure}

  Figure 2 gives the graphs of the AUMC, ALMC, and the lower bound $FT$ given by the Tverberg permanent conjecture for $r=6$.
  The value $h_3:=\max_{p\in[0,1]} h_{\Z^3}(p)$ is the $3$-dimensional monomer-dimer entropy of $\Z^3$ computed in \cite{FP}.
  \begin{figure}[here]
    \begin{center}
    \includegraphics[width=0.7\textwidth]{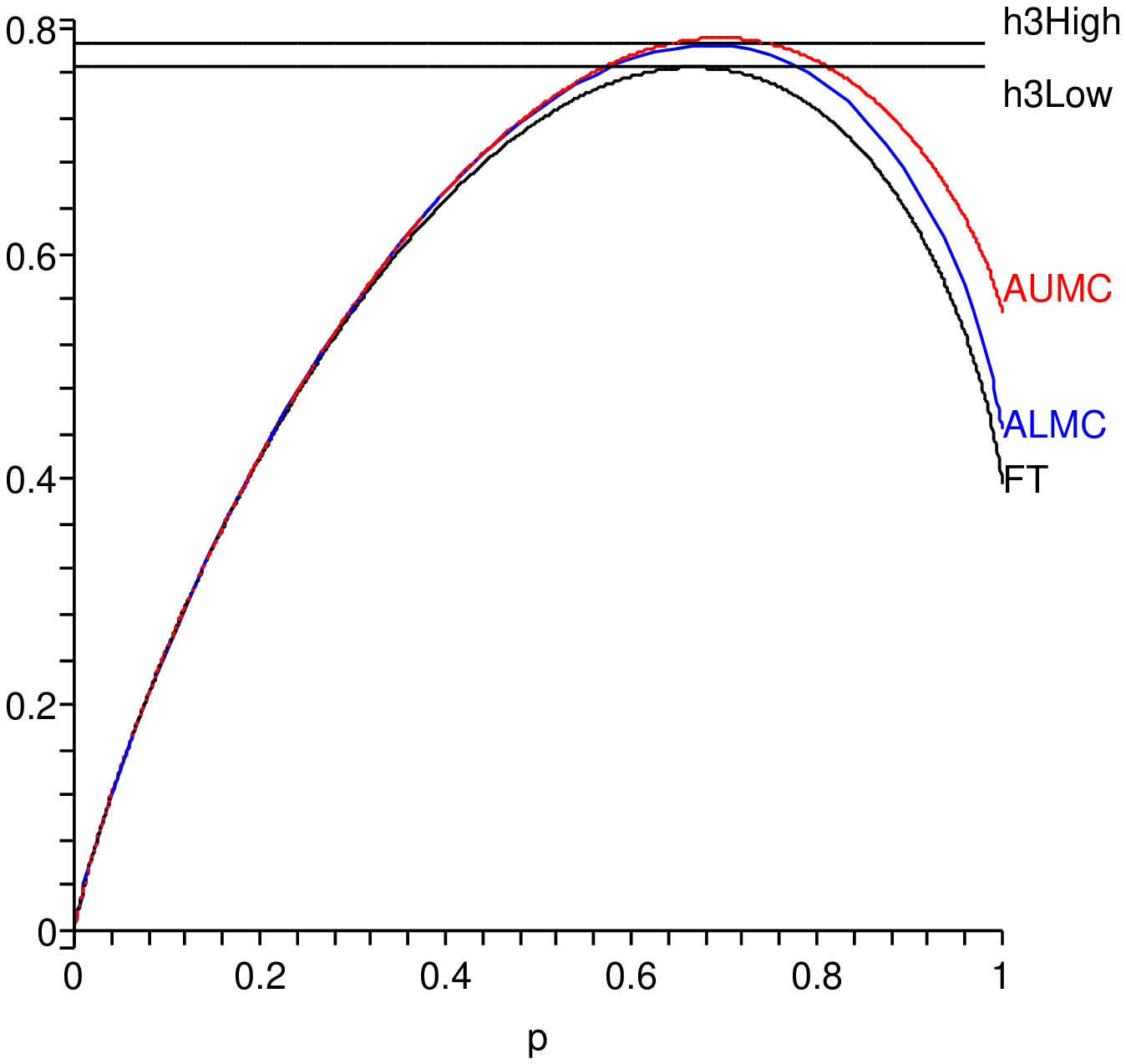}
    \caption{Monomer-dimer tiling of the $3$-dimensional grid:
    entropy as a function of dimer density. FT is the
    Friedland-Tverberg lower bound, h3Low and h3High are the known
    bounds for the monomer-dimer entropy. ALMC is the Asymptotic
    Lower Matching Conjecture. AUMC is the entropy of a countable
    union of $K_{6,6}$, conjectured to be an upper bound by the
    Asymptotic Upper Matching Conjecture.}\label{Fig:h3bounds}
    \end{center}
  \end{figure}

  Note that $h_l$  are very close to $\max_{p\in[0,1]} f(p,2l)(p)$ for $l=2,3$, where $f(p,2l)$ is given by \eqref{deffrp}, which is the graph
  corresponding to the ALMC.  The inequalities \eqref{lbfg08} yield the proof of the ALMC for $p=\frac{r}{r+s}$ for $s\in\N$ and $r\ge 3$
  \cite{FG08}.  It is shown in \cite{FF11}  that the ALMC holds for $h_{\Z^d}(p)$ for $d\ge 2$ and $p\in [0,1]$.  The proof of the ALMC is given
  in \cite{Gur11}.

  Figure 3 gives the continuous version of \eqref{uppkmatchbi}.  The green graph corresponds to the AUMC, i.e. the entropy of countable union of
  $K_{4,4}$.  The blue line, the graph of $\upp_{4,1}$ corresponds to the continuous version of ${n \choose k}^2(r!)^{\frac{k}{r}}$ in \eqref{uppkmatchbi}.  The orange line, the graph off $\upp_{4,2}$ corresponds to the continuous version of ${n \choose k} r^k$ in \eqref{uppkmatchbi}.  Note that $[0,1]$ is divided in two regions, in for which one the upper bounds $\upp_{4,j}$ is better than the other one.

    \begin{figure}[here]
      \includegraphics[width=0.7\textwidth]{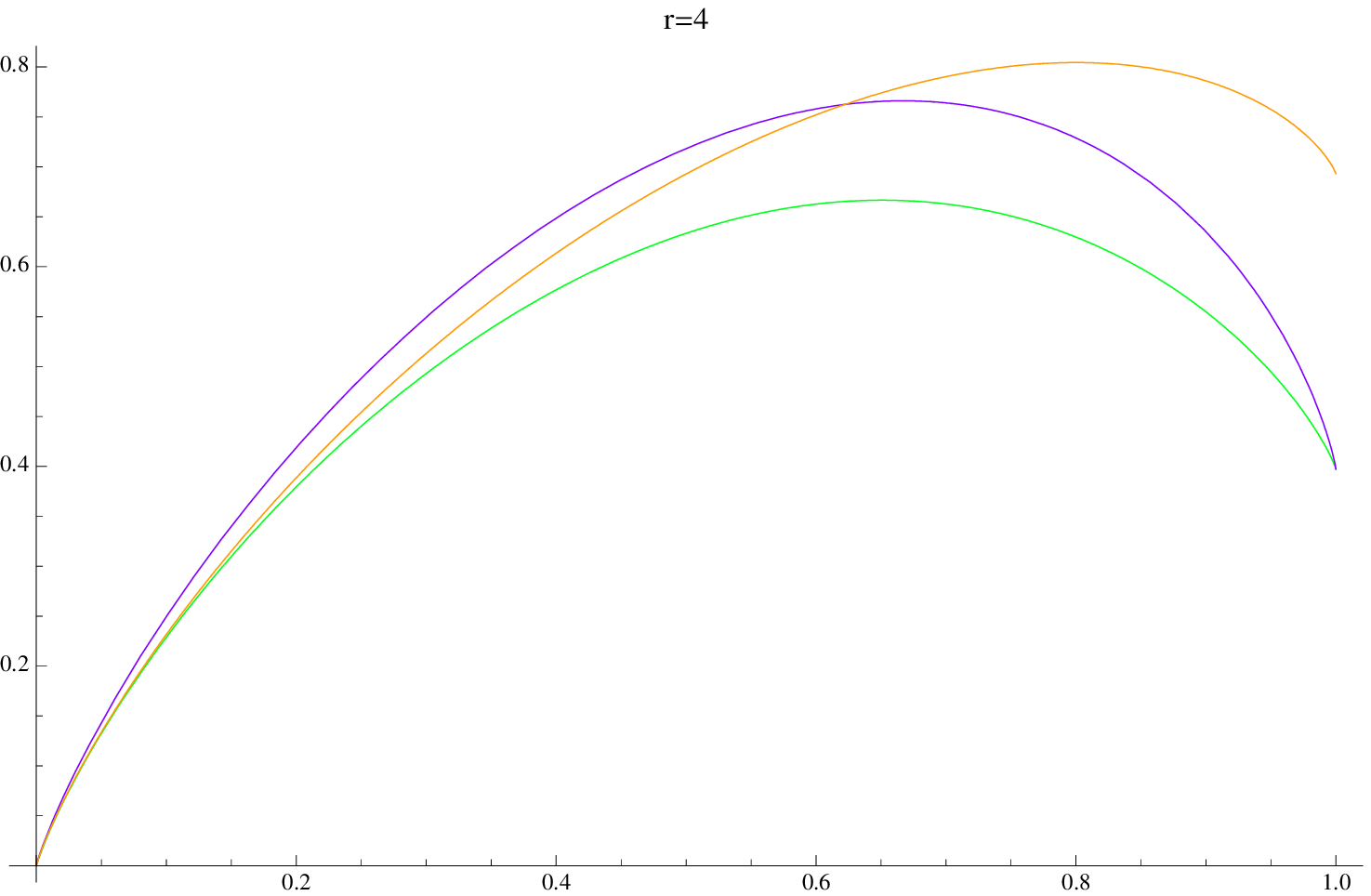}
      \caption{$h_{K(4)}$-green, $\upp_{4,1}$-blue, $\upp_{4,2}$-orange}
      \label{fig:f11}
 \end{figure}

 We close this section with open problems for infinite regular non-bipartite graphs.
 First one, find an analog of the LAMC for any infinite $r$-regular infinite graph for $p\in [0,\frac{2}{3}]$.
 Second, what is the analog of the ALMC for infinite cubic connected graphs with no bridges?

\end{document}